
  %

  \input amssym


  \font \bbfive		= bbm5
  \font \bbseven	= bbm7
  \font \bbten		= bbm10
  \font \eightbf	= cmbx8
  \font \eighti		= cmmi8 \skewchar \eighti = '177
  \font \eightit	= cmti8
  \font \eightrm	= cmr8
  \font \eightsl	= cmsl8
  \font \eightsy	= cmsy8 \skewchar \eightsy = '60
  \font \eighttt	= cmtt8 \hyphenchar\eighttt = -1

  \font \sixi		= cmmi6 \skewchar \sixi = '177
  \font \sixrm		= cmr6
  \font \sixsy		= cmsy6 \skewchar \sixsy = '60
  \font \tensc		= cmcsc10

  \font \titlefont	= cmbx12
  \scriptfont \bffam	= \bbseven
  \scriptscriptfont \bffam = \bbfive
  \textfont \bffam	= \bbten

  \newskip \ttglue

  \def \eightpoint {\def \rm {\fam 0 \eightrm }\relax
  \textfont 0= \eightrm
  \scriptfont 0 = \sixrm \scriptscriptfont 0 = \fiverm
  \textfont 1 = \eighti
  \scriptfont 1 = \sixi \scriptscriptfont 1 = \fivei
  \textfont 2 = \eightsy
  \scriptfont 2 = \sixsy \scriptscriptfont 2 = \fivesy
  \textfont 3 = \tenex
  \scriptfont 3 = \tenex \scriptscriptfont 3 = \tenex
  \def \it {\fam \itfam \eightit }\relax
  \textfont \itfam = \eightit
  \def \sl {\fam \slfam \eightsl }\relax
  \textfont \slfam = \eightsl
  \def \bf {\fam \bffam \eightbf }\relax
  \textfont \bffam = \bbseven
  \scriptfont \bffam = \bbfive
  \scriptscriptfont \bffam = \bbfive
  \def \tt {\fam \ttfam \eighttt }\relax
  \textfont \ttfam = \eighttt
  \tt \ttglue = .5em plus.25em minus.15em
  \normalbaselineskip = 9pt
  \def \MF {{\manual opqr}\-{\manual stuq}}\relax
  \let \sc = \sixrm
  \let \big = \eightbig
  \setbox \strutbox = \hbox {\vrule height7pt depth2pt width0pt}\relax
  \normalbaselines \rm }


  \def \ifundef #1{\expandafter \ifx \csname #1\endcsname \relax }


  \newcount \secno \secno = 0
  \newcount \stno \stno = 0
  \newcount \eqcntr \eqcntr = 0

  \def \track #1#2#3{\ifundef {#1}\else \hbox {\sixrm [#2\string #3] }\fi }

  \def \advseqnumbering {\global \advance \stno by 1 \global \eqcntr =0}

  \def \current {\number \secno \ifnum \number \stno = 0 \else .\number \stno \fi }

  \def \laberr #1{\message {*** RELABEL CHECKED FALSE for #1 ***}
      RELABEL CHECKED FALSE FOR #1, EXITING.
      \end }

  \def \deflabel #1#2{%
    \ifundef {#1}%
      \global \expandafter
      \edef \csname #1\endcsname {#2}%
    \else
      \edef \deflabelaux {\expandafter \csname #1\endcsname }%
      \edef \deflabelbux {#2}%
      \ifx \deflabelaux \deflabelbux \else \laberr {#1=(\deflabelaux )=(\deflabelbux )} \fi
      \fi
    \track {showlabel}{*}{#1}}

  \def \eqmark #1 {\advseqnumbering
    \eqno {(\current )}
    \deflabel {#1}{\current }}

  \def \subeqmark #1 {\global \advance \eqcntr by 1
    \edef \subeqmarkaux {\current .\number \eqcntr }
    \eqno {(\subeqmarkaux )}
    \deflabel {#1}{\subeqmarkaux }}

  \def \label #1 {\deflabel {#1}{\current }}
  \def \lcite #1{(#1\track {showlcit}{$\bullet $}{#1})}
  \def \forwardcite #1#2{\deflabel {#1}{#2}\lcite {#2}}


  \catcode `\@=11
  \def \c@itrk #1{{\bf #1}\track {showcitations}{\#}{#1}} 
  \def \c@ite #1{{\rm [\c@itrk{#1}]}}
  \def \sc@ite [#1]#2{[\c@itrk{#2}\hskip 0.7pt:\hskip 2pt #1]}
  \def \du@lcite {\if \pe@k [\expandafter \sc@ite \else \expandafter \c@ite \fi }
  \def \cite {\futurelet \pe@k \du@lcite }
  \catcode `\@=12


  \newcount \bibno \bibno = 0
  \newcount \bibtype \bibtype = 0 


  \def \newbib #1#2{\ifcase \number \bibtype
	\global \advance \bibno by 1 \edef #1{\number \bibno }\or
	\edef #1{#2}\or
	\edef #1{\string #1}\fi }

  \def \bibitem #1#2#3#4{\smallbreak \item {[#1]} #2, ``#3'', #4.}

  \def \references {\begingroup \bigbreak \eightpoint
    \centerline {\tensc References}
    \nobreak \medskip \frenchspacing }


  \def \Headlines #1#2{\nopagenumbers
    \headline {\ifnum \pageno = 1 \hfil
    \else \ifodd \pageno \tensc \hfil \lcase {#1} \hfil \folio
    \else \tensc \folio \hfil \lcase {#2} \hfil
    \fi \fi }}

  \long \def \Quote #1\endQuote {\begingroup \leftskip 35pt \rightskip 35pt
\parindent 17pt \eightpoint #1\par \endgroup }
  \long \def \Abstract #1\endAbstract {\bigskip \Quote \noindent #1\endQuote }

  \def \Note #1{\footnote {}{\eightpoint #1}}
  \def \Date #1 {\Note {\it Date: #1.}}


  \def \lcase #1{\edef \auxvar {\lowercase {#1}}\auxvar }

  \def \section #1 \par {\global \advance \secno by 1 \stno = 0
    \bigbreak \noindent {\bf \number \secno .\enspace #1.}
    \nobreak \medskip \noindent }

  \def \state #1 #2\par {\medbreak \noindent \advseqnumbering {\bf \current .\enspace #1.\enspace \sl #2\par }\medbreak }
  \def \definition #1\par {\state Definition \rm #1\par }

  \long \def \Proof #1\endProof {\medbreak \noindent {\it Proof.\enspace }#1
\ifmmode \eqno \endproofmarker $$ \else \hfill $\endproofmarker $ \looseness = -1 \fi \medbreak }

  \def \$#1{#1 $$$$ #1}
  \def \explain #1#2{\mathrel {\buildrel \hbox {\sixrm (#1)} \over #2}}
  \def \=#1{\explain {#1}{=}}

  \def \pilar #1{\vrule height #1 width 0pt}

  \newcount \footno \footno = 1
  \newcount \halffootno \footno = 1
  \def \footcntr {\global \advance \footno by 1
  \halffootno =\footno
  \divide \halffootno by 2
  $^{\number \halffootno }$}




  \def \Item #1{\smallskip \item {{\rm #1}}}
  \newcount \zitemno \zitemno = 0

  \def \izitem {\global \zitemno = 0}
  \def \zitemplus {\global \advance \zitemno by 1 \relax }
  \def \rzitem {\romannumeral \zitemno }
  \def \rzitemplus {\zitemplus \rzitem } 
  \def \zitem {\Item {{\rm (\rzitemplus )}}}
  \def \Zitem {\izitem \zitem }
  \def \zitemmark #1 {\deflabel {#1}{\rzitem }}

  \newcount \nitemno \nitemno = 0
  
  \def \nitem {\global \advance \nitemno by 1 \Item {{\rm (\number \nitemno )}}}

  \newcount \aitemno \aitemno = -1
  \def \boxlet #1{\hbox to 6.5pt{\hfill #1\hfill }}
  
  \def \aitemconv {\ifcase \aitemno a\or b\or c\or d\or e\or f\or g\or
h\or i\or j\or k\or l\or m\or n\or o\or p\or q\or r\or s\or t\or u\or
v\or w\or x\or y\or z\else zzz\fi }
  \def \aitem {\global \advance \aitemno by 1\Item {(\boxlet \aitemconv )}}
  \def \aitemmark #1 {\deflabel {#1}{\aitemconv }}


  \font \mf =cmex10
  \def \union {\mathop {\raise 9pt \hbox {\mf S}}\limits }
  \def \inters {\mathop {\raise 9pt \hbox {\mf T}}\limits }

  \def \<{\left \langle \vrule width 0pt depth 0pt height 8pt }
  \def \>{\right \rangle }
  \def \({\big (}
  \def \){\big )}
  
  \def \and {\hbox {,\quad and \quad }}
  \def \calcat #1{\,{\vrule height8pt depth4pt}_{\,#1}}

  \def \imply {\kern 7pt \Rightarrow \kern 7pt}
  \def \for #1{,\quad \forall \,#1}
  \def \endproofmarker {\square } 
  \def \"#1{{\it #1}\/} \def \umlaut #1{{\accent "7F #1}}
  
  \def \*{\otimes }
  \def \caldef #1{\global \expandafter \edef \csname #1\endcsname {{\cal #1}}}
  \def \bfdef #1{\global \expandafter \edef \csname #1\endcsname {{\bf #1}}}
  \bfdef N \bfdef Z \bfdef C \bfdef R


  \def \Caixa #1{\setbox 1=\hbox {$#1$\kern 1pt}\global \edef \tamcaixa {\the \wd 1}\box 1}
  \def \caixa #1{\hbox to \tamcaixa {$#1$\hfil }}



  \catcode `\@=11

  \def \overparenOnefill {$\m@th
  \setbox 0=\hbox {$\braceld $}%
  \braceld \leaders \vrule height\ht 0 depth0pt\hfill
  \leaders \vrule height\ht 0 depth0pt\hfill \bracerd $}

  \def \overparenOne #1{\mathop {\vbox {\m@th\ialign {##\crcr \noalign {\kern -1pt}
  \overparenOnefill \crcr \noalign {\kern 3pt\nointerlineskip }
  $\hfil \displaystyle {#1}\hfil $\crcr }}}\limits }

  \def \overparenTwofill {$\m@th
  \lower 0.3pt \hbox {$\braceld $}
  \leaders \vrule depth 0pt height1pt \hfill
  \lower 0.3pt \hbox {$\bracerd $}$}

  \def \overparenTwo #1{\mathop {\vbox {\ialign {##\crcr \noalign {\kern -1pt}
  \overparenTwofill \crcr \noalign {\kern 3pt\nointerlineskip }
  $\hfil \displaystyle {#1}\hfil $\crcr }}}\limits }

  \catcode `\@=12


\font\rs=rsfs10

\caldef G
\def \ScriptCapitalF{{\cal F}}  \def \ScriptCapitalF{\hbox{\rs F}}
\def \ScriptCapitalH{{\cal H}}
\def \CGx{C^*\big(\G(x)\big)}
\def \gt#1{\widehat{#1}}
\def \zero{^{(0)}}
\def \Gz{
  \G\zero}

\def \bool#1{[{\scriptstyle #1}]}
\def \soma{\mathop {\textstyle\sum}\limits}

\def \.{\odot}
\def \etale{\'etale}
\def \L{L}
\def \HL{H_\L}
\def \Hi#1{H_{{\rm Ind}\kern0.5pt #1}}
\def \sysind#1#2#3{{\rm Ind}_{\scriptscriptstyle #1}^{\scriptscriptstyle #2}\kern0.5pt #3}
\def \ind#1{\sysind HG#1}
\def \fix{\noindent $\blacktriangleright $\quad}
\def \s{s} \def \r{r} \def \e#1{e_{#1}}
\def \star{{*\kern 0.7pt}}
\def \Hx{H_x}
\def \X{M}
\def \Ker{\hbox{Ker}}
\def \ustar{\!_{_*}}

 \bibtype 0

\newbib{\Allan}{A}
\newbib {\ClaireRenault}{AR}
\newbib {\Amelia}{BFK}
\newbib {\Dougkas}{D}
\newbib {\IW}{IW}
\newbib {\IWTwo}{IW2}
\newbib {\KarlovOne}{K}
\newbib {\MRW}{MRW}
\newbib {\renault}{Re}
\newbib {\Roch}{Ro}

  \Headlines {Invertibility in groupoid C*-algebras}
    {R.~Exel}
  \null\vskip -1cm
  \centerline{\titlefont INVERTIBILITY IN GROUPOID C*-ALGEBRAS}

  \bigskip
  \centerline
  {\tensc R.~Exel\footnote{*}{\eightrm Partially supported by CNPq.}}


  \footnote{\null}
  {\eightrm 2010 \eightsl Mathematics Subject Classification:
  \eightrm
  22A22, 
  46L05, 
  46L55. 
  }

  \bigskip

  \Abstract
  Given a second-countable, Hausdorff, {\etale}, amenable groupoid $\G$ with compact unit space, we show that an element
$a$ in $C^*(\G)$ is invertible if and only if $\lambda_x(a)$ is invertible for every $x$ in the unit space of $\G$, where $\lambda_x$
refers to the \"{regular representation} of $C^*(\G)$ on $\ell_2(\G_x)$.  We also prove that, for
every $a$ in $C^*(\G)$, there exists some $x\in\Gz$ such that  $\Vert a\Vert = \Vert\lambda_x(a)\Vert.$
  \endAbstract

  \section Introduction

The structure of certain C*-algebras is often best studied via large families of *-represen\-ta\-tions.  According to this
point of view, one tries to deduce the properties of any given element of the algebra by means of the properties of its
images under the representations provided.
  Here we shall mostly be interested in
\"{invertibility} questions, and thus on families of representations of a given C*-algebras which are large enough to
determine when an element is invertible.

One of the first, and arguably also the  most influential such result is the  Allan-Douglas
local principle  \cite[Corollary  2.10]{\Allan},
\cite[Theorem 7.47]{\Dougkas}, which asserts that an element in a unital Banach algebra is invertible if and only if it
is invertible modulo certain ideals associated to the points of the spectrum of a given central subalgebra.
This principle has been generalized to \"{nonlocal} algebras   (see \cite{\KarlovOne} and the references given there) and
has successfully been applied to study Fredholm singular integral operators with semi-almost periodic coefficients
\cite{\Amelia}.

The present paper is an attempt to transpose  the \"{local-trajectory method} of \cite{\KarlovOne} to the context of
groupoid C*-algebras.  Since invertibility only makes sense on unital algebras, and since
the C*-algebra of a groupoid is unital only when the groupoid is \'etale and has a compact
unit space, we restrict ourselves to this case (however our work suggests questions that might be relevant for more
general groupoids).
To be precise, our main result, Theorem \forwardcite {MainThm}{2.10},
applies to  second-countable, Hausdorff, {\etale}, amenable groupoids with compact unit
space.  Given such a groupoid $\G$, we show that an element $a$ in the groupoid C*-algebra $C^*(\G)$ is invertible if and
only, for every $x$ in the unit space of $\G$, one has that $\lambda_x(a)$ is invertible,  where $\lambda_x$
is the \"{regular representation} of $C^*(\G)$ on $\ell_2(\G_x)$.

A crucial tool used to prove our main result is the theory of induced representations started by Renault in
\cite[Chap. II, \S2]{\renault} and improved by Ionescu and Williams in \cite{\IW} and \cite{\IWTwo}.

Recall that the amenability assumption on $\G$ implies \cite[Theorem 6.1.4.(iii)]{\ClaireRenault} that
  $$
  \Vert a\Vert = \sup_{x\in\Gz}\Vert\lambda_x(a)\Vert   \for a\in C^*(\G).
  \eqmark IntrNorming
  $$
  As a byproduct of our work  we have found  a small improvement of this result, namely
Corollary \forwardcite {MainCorol}{3.4}, below, which asserts that
  $$
  \Vert a\Vert = \max_{x\in\Gz}\Vert\lambda_x(a)\Vert  \for a\in C^*(\G),
  \eqmark StrictNorming
  $$
  which is to say that
the  supremum in \lcite{\IntrNorming} is in fact \"{attained} for every $a$.  The proof of this fact is a
straightforward  combination of Theorem
\forwardcite {MainThm}{2.10} with a result of S. Roch \cite{\Roch}, which we carefully describe below.

Even though the invertibility question treated in \forwardcite {MainThm}{2.10} only makes sense for groupoids with
compact  unit space, \lcite{\StrictNorming} applies to a wider context.  A sensible question to be asked at this point is therefore whether
or not \lcite{\StrictNorming} holds in the absence of the compactness hypothesis.

Dropping the assumption that $\G$ is amenable, it is well known that \lcite{\IntrNorming} holds as long as we replace
the full by the reduced groupoid C*-algebra.  So it makes sense to ask whether or not
  $$
  \Vert a\Vert = \max_{x\in\Gz}\Vert\lambda_x(a)\Vert  \for a\in C^*_r(\G)\ ?
  \eqmark StrictNormingNONAM
  $$
  Unfortunately we have not been able to answer any of these questions, which we are then forced to  leave as open problems.

Attaining the supremum is a well known property of continuous functions on compact spaces, so a proof of
\lcite{\StrictNormingNONAM} could be obtained, at least in the case of a compact unit space, should we be able to prove
that the function
  $$
  x \mapsto \Vert\lambda_x(a)\Vert
  $$
  is continuous for every $a\in C^*_r(\G)$.  However sensible this appears to be, we have not been able to determine its
validity.

Last, but not least, I would like to thank Am\'elia Bastos and the members of ``The Center for Functional Analysis and
Applications - CEAF'' of the ``Instituto Superior T\'ecnico de Lisboa'' for bringing their work to my attention and also
for their warm hospitality during two visits there where many interesting conversations on these topics took place and
where the ideas for the present work developed.  I would also like to thank Jean Renault for helpful e-mail exchanges.

  \section Sufficient family of representations

Let $A$ be  a unital C*-algebra.  The following concept appears in \cite[Section 5]{\Roch}.

\definition A family $\ScriptCapitalF$ of non-degenerated representations (always assumed to preserve the involution) of $A$ is called
\"{sufficient} if, for every $a$ in $A$, one has that
  $$
  a \hbox{ is invertible} \iff \pi(a) \hbox{ is invertible for all } \pi\in\ScriptCapitalF.
  $$

Observe that the implication ``$\Rightarrow$" is always true, so the relevant property conveyed by this definition is the
implication ``$\Leftarrow$".

\state Proposition \label Invertibility The set of all irreducible representations of $A$ is a sufficient family of
representations.

\Proof If $a$ is a non-invertible element of $A$, then either $a^*a$ or $aa^*$ are non-invertible.  So we may assume,
without loss of generality that $a^*a$ is non-invertible.  Let $B$ be the closed *-subalgebra of $A$ generated by $a^*a$
and $1$, and let $X$ be the compact spectrum of $B$.  Since $a^*a$ is non-invertible, there exists some point $x_0$ in $X$
such that $\gt{a^*a}(x_0) = 0$, where the hat indicates the Gelfand transform.

The map
  $$
  \phi: b \in B \mapsto\gt b(x_0) \in{\bf C}
  $$
  is therefore a pure state of $B$, which may be extended to a pure state $\psi$ on $A$.  Let $\pi$ be the GNS representation
associated to $\psi$, so that $\pi$ is an irreducible representation.  If $\xi$ is the associated cyclic vector we have
  $$
  \Vert\pi(a)\xi\Vert^2 = \langle\pi(a)\xi,\pi(a)\xi\rangle= \langle\pi(a^*a)\xi,\xi\rangle= \psi(a^*a) =
  \phi(a^*a) = \gt{a^*a}(x_0) = 0.
  $$
  It follows that the operator $\pi(a)$ is not injective and hence non-invertible.
  \endProof


\fix From now on we will be interested in the question of sufficiency for groupoid C*-algebras.   We therefore fix a
second-countable, Hausdorff, {\etale} groupoid  $\G$, with source and range maps denoted by ``$s$'' and ``$r$'', respectively.

\bigskip Given $x$ in the unit space $\Gz$ of $\G$, we shall use the following standard notations:
  $$
  \def \crr{\hfill\cr\pilar{18pt}}
  \matrix{
  \G_x & = & \{\gamma\in\G: \s(\gamma) = x\}, \crr
  \G^x & = & \{\gamma\in\G: \r(\gamma) = x\}, \quad \hbox{and}\crr
  \G(x) & = & \G_x \cap\G^x.\hfill
  }
  $$

Consider the Hilbert space $\Hx = \ell_2(\G_x)$ and the \"{regular representation} $\lambda_x$ of $C_c(\G)$ on $H_x$, given by
  $$
  \lambda_x(f)\xi\calcat \gamma= \soma_{\gamma'\gamma'' = \gamma} f(\gamma')\xi(\gamma'')
  \for f \in C_c(\G) \for \xi\in\Hx \for\gamma\in\G_x,
  $$
  which is well known to extend to $C^*(\G)$.  For each $\gamma$ in $ \G_x$, let $\e\gamma$ be the basis vector of $H_x$ corresponding to $\gamma$.

\state Proposition \label MatrixEntries For every $\gamma_1$ and $\gamma_2$ in $\G_x$, and all $f$ in $C_c(\G)$, one has that
  $$
  \langle\lambda_x(f)\e{\gamma_1},\e{\gamma_2}\rangle= f(\gamma_2\gamma_1^{-1}).
  $$

\Proof
  We have
  $$
  \langle\lambda_x(f)\e{\gamma_1},\e{\gamma_2}\rangle= \lambda_x(f)\e{\gamma_1}\calcat{\gamma_2} =
  \soma_{\gamma'\gamma'' = \gamma_2} f(\gamma')\e{\gamma_1}(\gamma'') = \soma_{\gamma'\gamma_1 = \gamma_2} f(\gamma') = f(\gamma_2\gamma_1^{-1}).  \endProof

\state Proposition
  Let $\ScriptCapitalH$ be a closed sub-groupoid of $\G$, viewed as a topological groupoid with the relative topology.  Then the
following are equivalent:
  \Zitem the restriction of the range map $r$ to $\ScriptCapitalH$, viewed as a mapping
  $$
  r|_\ScriptCapitalH: \ScriptCapitalH \to\ScriptCapitalH\zero,
  $$
  is an open mapping,
  \zitem $\ScriptCapitalH$ is \etale.

  \Proof Assuming (i),
  let $\gamma\in\ScriptCapitalH$ and choose an open set $U \subseteq\G$ such that $r$ is a homeomorphism from $U$ onto the open set $r(U) \subseteq\Gz.$
Then $U\cap\ScriptCapitalH$ is open in the relative topology of $\ScriptCapitalH$ and, by (i), we have that $r(U\cap\ScriptCapitalH)$ is open in $\ScriptCapitalH\zero$.  It is
then clear that $r$ is a homeomorphism from $U\cap\ScriptCapitalH$ to $r(U\cap\ScriptCapitalH)$, showing that $r|_\ScriptCapitalH$ is a local homeomorphism and hence
that $\ScriptCapitalH$ is \etale.  The converse is evident.
  \endProof

\fix From now on we fix a closed sub-groupoid $\ScriptCapitalH \subseteq\G$, satisfying the equivalent conditions above.  We will denote the
unit spaces of $\G$ and $\ScriptCapitalH$ as follows
  $$
  X: = \Gz
  \and
  Y: = \ScriptCapitalH\zero.
  $$
  Since $H$ is closed in $\G$ and since $Y = \ScriptCapitalH\cap X$, we see that $Y$ is a closed subspace of $X$.

Let us briefly describe the process of inducing representations from $C^*(\ScriptCapitalH)$ to $C^*(\G)$,
  cf.~\cite[Chap. II, \S2]{\renault} and
  \cite[Section 2]{\IWTwo}.
  Given a representation $\L$ of $C^*(\ScriptCapitalH)$ on a Hilbert space $\HL$, we want to produce a representation $\ind L$ of
$C^*(\G)$ on a Hilbert space $\Hi L$.  In order to do so, consider the closed subset of $\G$ given by
  $$
  \G_Y = s^{-1}(Y) = \{\gamma\in\G: \s(\gamma) \in Y\}.
  $$
  For $\varphi$ and $\psi$ in $C_c(\G_Y)$, define $\langle\varphi,\psi\rangle\ustar$ in $C_c(\ScriptCapitalH)$, by
  $$
  \langle\varphi,\psi\rangle\ustar(\zeta) = \sum_{\gamma_1\gamma_2 = \zeta}\overline{\varphi(\gamma_1^{-1})}\psi(\gamma_2)
  \for \zeta\in\ScriptCapitalH.
  $$
  It should be noticed that the above sum ranges over all pairs of elements $\gamma_1$ and $\gamma_2$ in $\G$ (as opposed to $\ScriptCapitalH$),
whose product equals $\zeta$.  In this case notice that both $\r(\gamma_1)$ and $\s(\gamma_2)$ lie in $Y$, so that $\gamma_1^{-1}$ and $\gamma_2$ indeed
belong to the domain of $\varphi$ and $\psi$, respectively.

By \cite[Theorem 2.8]{\MRW}, one has that in fact $C_c(\G_Y)$ may be completed to a right $C^*(\ScriptCapitalH)$--Hilbert module, which
we will denote by $\X$, the appropriate right multiplication being that which is described in \cite[page 11]{\MRW}.  It
is therefore profitable to view $\langle\cdot\,,\cdot\rangle\ustar$ as a $C^*(\ScriptCapitalH)$--valued map.

The space $\Hi L$, on which the induced representation will act, is then defined to be the completion of
  $$
  C_c(\G_Y) \otimes\HL,
  $$
  relative to the inner-product
  $$
  \langle\varphi\otimes\xi, \psi\otimes\eta\rangle: = \big\langle L\big(\langle\psi,\varphi\rangle\ustar\big) \xi,\eta\big\rangle
  \for \varphi, \psi\in C_c(\G_Y) \for \xi, \eta\in\HL.
  $$
  One next gives $C_c(\G_Y)$ the structure of a left $C_c(\G)$--module by setting
  $$
  (f\star\varphi)(\gamma) : = \sum_{\gamma_1\gamma_2 = \gamma}f(\gamma_1)\varphi(\gamma_2)
  \for f \in C_c(\G) \for \varphi\in C_c(\G_Y)\for \gamma\in\G_Y.
  $$
  Again by \cite[Theorem 2.8]{\MRW}, the above left-module structure may be extended to a bounded multiplication
operation
  $$
  (a,x) \in C^*(\G) \times\X \mapsto ax \in\X.
  $$

In order to define the induced representation one may either work with the completion $\X$ described above or take the
more pedestrian point of view of sticking to compactly supported functions.  Taking the latter approach,
  for $f \in C_c(\G)$ one initially defines
  $\ind L(f)$ on the dense subspace
  $C_c(\G_Y) \otimes\HL \subseteq\Hi L$, by the formula
  $$
  \ind L(f)(\varphi\otimes\xi) : = (f\star\varphi) \otimes\xi
  \for \varphi\in C_c(\G_Y) \for \xi\in\HL,
  $$
  and then extend it by continuity to $\Hi L$.  This provides a *-representation of $C_c(\G)$ on $\Hi L$ which, in turn,
may be extended to the whole of $C^*(\G)$.

The resulting representation of $C^*(\G)$ on $\Hi L$ is denoted by $\ind L$, and is called the \"{representation induced
by $L$ from $\ScriptCapitalH$ up to $\G$}.  For more details, see
  \cite[Chap. II, \S2]{\renault} and
  \cite[Section 2]{\IWTwo}.


\bigskip\fix Fix, for the time being, an element $x \in X$.

\bigskip We would now like to consider the question of inducing representations from $\ScriptCapitalH: = \G(x)$ up to $\G$.
Observing that
  $$
  Y = \G(x)\zero = \{x\},
  $$
  we have that $\G_Y = \G_x$, which is a discrete topological space.  Consequently $C_c(\G_Y)$ is linearly generated by
the set
  $$
  \{\e \gamma: \gamma\in\G_x\},
  $$
  where $\e\gamma$ denotes the characteristic function of the singleton $\{\gamma\}$.

\state Proposition \label basisForSpace Given $\gamma,\gamma' \in\G_x$, we have that
  $$
  \langle\e\gamma, \e{\gamma'}\rangle\ustar =
  \left\{\matrix{
  \delta_{\gamma^{-1}\gamma'}, & \hbox{ if } \r(\gamma) = \r(\gamma'), \cr \pilar{15pt}
  \hfill 0 \hfill, & \hbox { otherwise,}\hfill
  }\right.
  $$
  where, for each $h \in\G(x)$, we denote by $\delta_h$ the characteristic function of the singleton $\{h\}$, viewed as an
element of $C_c\big(\G(x)\big) \subseteq \CGx$.

\Proof We have, for every $\zeta\in\G(x)$, that
  $$
  \langle\e\gamma, \e{\gamma'}\rangle\ustar(\zeta) =
  \sum_{\gamma_1\gamma_2 = \zeta}\overline{\e\gamma(\gamma_1^{-1})}\e{\gamma'}(\gamma_2) =
  \bool{\gamma^{-1}\gamma' = \zeta},
  $$
  where the brackets denote the Boolean value of the expression inside, with the convention that a syntactically
incorrect expression, e.g.~when the multiplication $\gamma^{-1}\gamma'$ is illegal, the value is zero.

  Thus, when $\r(\gamma) = \r(\gamma')$, we have that the product $\gamma^{-1}\gamma'$ is defined, evidently giving an element of $\G(x)$ and,
in this case,
  $$
  \langle\e\gamma, \e{\gamma'}\rangle\ustar = \delta_{\gamma^{-1}\gamma'},
  $$
  On the other hand, when $\r(\gamma) \neq\r(\gamma')$, we clearly have that
  $
  \langle\e\gamma, \e{\gamma'}\rangle\ustar = 0.
  $
  \endProof

The following elementary result is included in order to illustrate a simple example.

\state Proposition \label ExampleInduceRegular \rm
  Let $\Lambda$ be the left-regular representation of $\CGx$ on $\ell_2(\G(x))$.  Then $\sysind{\G(x)}{\G}\Lambda$ is unitarily
equivalent to $\lambda_x$.

\Proof For each element $\gamma\in\G_x$, and each $g \in\G(x)$, consider the element
  $$
  \varphi_{\gamma,g} = \e\gamma\otimes\e g \in C_c(\G_x)\otimes\ell_2(\G(x)) \subseteq\Hi \Lambda.
  $$
  We first claim that
  $$
  \langle\varphi_{\gamma,g},\varphi_{\gamma',g'}\rangle= \bool{\gamma g = \gamma'g'} \for \gamma,\gamma' \in\G_x \for g,g' \in\G(x).
  \subeqmark InnerTensorProd
  $$
  In fact, we have
  $$
  \langle\varphi_{\gamma,g},\varphi_{\gamma',g'}\rangle= \big\langle\e\gamma\otimes\e{\r(\gamma)},\e{\gamma'}\otimes\e{\r(\gamma')}\big\rangle=
  \big\langle\Lambda\big(\langle\e{\gamma'},\e\gamma\rangle\ustar\big)\e{\r(\gamma)},\e{\r(\gamma')}\big\rangle= (\dagger)
  $$
  Consequently, when $\r(\gamma) \neq\r(\gamma')$ we have by \lcite{\basisForSpace} that
  $\langle\varphi_{\gamma,g},\varphi_{\gamma',g'}\rangle= 0$, which proves \lcite{\InnerTensorProd} in this case.
  If $\r(\gamma) = \r(\gamma')$ then, again by \lcite{\basisForSpace}, it follows that
  $$
  (\dagger) = \big\langle\Lambda\big(\delta_{{\gamma'}^{-1}\gamma}\big)\e g,\e{g'}\big\rangle=
  \big\langle\e{{\gamma'}^{-1}\gamma g},\e{g'}\big\rangle=
  \bool{{\gamma'}^{-1}\gamma g = g'} = \bool{\gamma g = \gamma'g'},
  $$
  proving \lcite{\InnerTensorProd}.
  In particular, this implies that
  $$
  \langle\varphi_{\gamma,g},\varphi_{\gamma',g'}\rangle= \langle\varphi_{\gamma g,x},\varphi_{\gamma',g'}\rangle,
  $$
  and since the collection of all $\varphi_{\gamma',g'}$ evidently spans $\Hi \Lambda$, we have that $\varphi_{\gamma,g} = \varphi_{\gamma g,x}$, and it is then
clear that the mapping
  $$
  \e\gamma\mapsto\varphi_{\gamma,x}
  $$
  extends to a unitary operator
  $U: \Hx \to\Hi \Lambda$.
  Given $f \in C_c(\G)$, we claim that
  $$
  \big\langle U^*(\ind \Lambda(f))U\e\gamma,\,\e{\gamma'}\big\rangle= \langle\lambda_x(f)\e\gamma,\e{\gamma'}\rangle
  \for\gamma,\gamma' \in\G_x.
  \subeqmark ToProveEquivalent
  $$
  In order to verify it observe that the left-hand side equals
  $$
  \big\langle\ind \Lambda(f)(\varphi_{\gamma,x}),\,\varphi_{\gamma',x}\big\rangle=
  \big\langle(f\star\e\gamma)\otimes\e x,\,\e{\gamma'}\otimes\e {x'}\big\rangle=
  \big\langle\Lambda\big(\langle\e{\gamma'},f\star\e\gamma\rangle\ustar\big)\e x,\,\e x\big\rangle= (\diamondsuit)
  $$
  After checking that
  $$
  f\star\e\gamma= \sum_{\eta\in\G_x}f(\eta\gamma^{-1})\e\eta,
  $$
  we conclude that
  $$
  (\diamondsuit) =
  \sum_{\eta\in\G_x} f(\eta\gamma^{-1}) \big\langle\Lambda\big(\langle\e{\gamma'},\e\eta\rangle\ustar\big)\e x,\,\e x\big\rangle=
  \sum_
    {\scriptstyle \eta\in\G_x \atop \scriptstyle \r(\gamma') = \r(\eta)}
    f(\eta\gamma^{-1}) \big\langle\Lambda\big(\gamma_{{\gamma'}^{-1}\eta}\big)\e x,\,\e x\big\rangle\$ =
  \sum_
    {\scriptstyle \eta\in\G_x \atop \scriptstyle \r(\gamma') = \r(\eta)}
    f(\eta\gamma^{-1}) \big\langle\e{{\gamma'}^{-1}\eta},\,\e x\big\rangle=
  f(\gamma'\gamma^{-1}) \ = {\MatrixEntries}
  \langle\lambda_x(f)\e\gamma,\e{\gamma'}\rangle.
  $$
  This proves \lcite{\ToProveEquivalent}, and since $\gamma$ and $\gamma'$ are arbitrary, we conclude that $U^*(\ind \Lambda(f))U =
\lambda_x(f)$, finishing the
  proof.
  \endProof

Notice that there are two completions of $C_c(\G_x)$ which are relevant to us.  On the one hand $\X$ is the completion
under the $\CGx$--valued inner-product $\langle\cdot\,,\cdot\rangle\ustar$, and, on the other, $\Hx$ is the completion for the 2-norm.  These
two spaces are related to each other by the following.

\state Proposition
  There is a bounded linear map
  $$
  j:\X \to\Hx,
  $$
  such that $j(\varphi) = \varphi$, for every $\varphi\in C_c(\G_x)$.

\Proof
  Given $\varphi\in C_c(\G_x)$, notice that
  $$
  \Vert\varphi\Vert_2^2 = \sum_{\gamma\in\G_x}\overline{\varphi(\gamma)}\varphi(\gamma) = \langle\varphi,\varphi\rangle\ustar(1) \leq
  \Vert\langle\varphi,\varphi\rangle\ustar\Vert_{_{C^*(\G(x))}} = \Vert\varphi\Vert_\X^2.
  $$
  This implies that the identity map on $C_c(\G_x)$ is continuous for $\Vert{\cdot}\Vert_\X$ on its domain and the 2-norm on its
codomain.  The required map is then obtained by a continuous extension.  \endProof

If $\zeta\in\G(x)$, we have a well defined bijective map
  $$
  \gamma\in\G_x \mapsto\gamma\zeta\in\G_x,
  $$
  and hence the map
  $$
  R_\zeta:\Hx \to\Hx ,
  $$
  defined by
  $$
  R_\zeta(\xi)\calcat \gamma= \xi(\gamma\zeta)
  \for \xi\in\Hx \for \gamma\in\G_x,
  $$
  is a unitary operator.  It is also easy to see that $R_{\zeta_1}\!\circ R_{\zeta_2} = R_{\zeta_1\zeta_2}$, which is to say that $R$ is a
unitary representation of $\G(x)$ on $\Hx $.

This representation will play an important role in our next result, but before stating it, we need to introduce a
notation.

Given any discrete group $G$, and any $\zeta\in G$, the map
  $$
  f \in C_c(G) \mapsto f(\zeta) \in\C
  $$
  is well known to extend to a bounded linear functional on $C^*(G)$, which we will denote by
  $$
  a \in C^*(G) \mapsto\hat a(\zeta) \in\C.
  $$

\state Proposition \label MainIdentity For every $a \in C^*(\G)$, every $x,y \in\X $, and every $\zeta\in\G(x)$, we have that
  $$
  \widehat {\langle x,a y\rangle\ustar}(\zeta) = \Big\langle\lambda_x(a)R_\zeta\big(j(y)\big),j(x)\Big\rangle.
  $$

\Proof Given $f \in C_c(\G)$, and $\psi,\varphi\in C_c(\G_x) $, we have
  $$
  \langle\varphi,f\star\psi\rangle\ustar(\zeta) =
  \sum_{\gamma_1\gamma_2 = \zeta}\overline{\varphi(\gamma_1^{-1})}(f\star\psi)(\gamma_2) =
  \sum_{\gamma_1\gamma_2\gamma_3 = \zeta}\overline{\varphi(\gamma_1^{-1})}f(\gamma_2)\psi(\gamma_3) = \cdots
  $$
  With the change of variables
  ``$
  \gamma_3' = \gamma_3 \zeta^{-1}
  $"
  the above equals
  $$
  \cdots=
  \sum_{\gamma_1\gamma_2\gamma_3' = x}\overline{\varphi(\gamma_1^{-1})}f(\gamma_2)\psi(\gamma_3'\zeta) =
  \sum_{\gamma_1\gamma_2\gamma_3' = x}\overline{\varphi(\gamma_1^{-1})}f(\gamma_2)R_\zeta(\psi)(\gamma_3') =
  \big\langle f\star R_\zeta(\psi),\varphi\big\rangle.
  $$
  This gives that
  $$
  \langle\varphi,f\star\psi\rangle\ustar(\zeta) = \big\langle f\star R_\zeta(\psi),\varphi\big\rangle,
  $$
  and the proof is concluded upon replacing
  \medskip \item{$\bullet$} $f$ by the terms of a sequence $\{f_n\}_n$ converging to $a$ in $\CGx$,
  \medskip \item{$\bullet$} $\varphi$ by the terms of a sequence $\{\varphi_n\}_n$ converging to $x$ in $\X$, and finally
  \medskip \item{$\bullet$} $\psi$ by the terms of a sequence $\{\psi_n\}_n$ converging to $y$ in $\X$.
  \endProof

\medskip \state Corollary \label Comparison
  Given $x \in X$, suppose that $a$ is an element of $C^*(\G)$ such that $\lambda_x(a) = 0$.  Then
  $$
  \sysind{\G(x)}{\G}L(a) = 0,
  $$
  for any representation $L$ of $\CGx$ which is weakly contained in $\Lambda$.

\Proof By \lcite{\MainIdentity}, we deduce that
  $$
  \widehat {\langle x,a y\rangle\ustar}(\zeta) = 0
  \for \zeta\in\G(x)
  \for x,y \in\X.
  $$
  Temporarily fixing $x$ and $y$, we then deduce that
  $\Lambda\big(\langle x,a y\rangle\ustar\big) = 0$, and hence that
  $$
  L\big(\langle x,a y\rangle\ustar\big) = 0, \subeqmark Zerou
  $$
  for any $L$ as in the statement.
  Given $f \in C_c(\G)$, $\varphi,\psi\in C_c(\G_x)$ and $\xi,\eta\in\HL$, we have that
  $$
  \big\langle\sysind{\G(x)}{\G}L(f) (\varphi\otimes\xi), \psi\otimes\eta\big\rangle= \big\langle(f\star\varphi)\otimes\xi, \psi\otimes\eta\big\rangle= \big\langle L\big(\langle\psi,f\star\varphi\rangle\ustar\big)\xi,\eta\big\rangle.
  $$
  Applying this for $f$ ranging in a sequence $\{f_n\}_n$ converging to $a$ in $\CGx$, we conclude that
  $$
  \big\langle\sysind{\G(x)}{\G}L(a) (\varphi\otimes\xi), \psi\otimes\eta\big\rangle= \big\langle L\big(\langle\psi,a\varphi\rangle\ustar\big)\xi,\eta\big\rangle\ = {\Zerou} 0,
  $$
  from where the conclusion follows easily.  \endProof

We may now prove our main result:

\state Theorem \label MainThm Let $\G$ be a second-countable, Hausdorff, {\etale} groupoid, such that $\Gz$ is compact.  Suppose
moreover that $\G$ is amenable.  Then $\{\lambda_x\}_{x\in\Gz}$ is a sufficient  family of representations for $C^*(\G)$.  In other
words, if $a \in C^*(\G)$ is such that $\lambda_x(a)$ is invertible for every $x$ in the unit space of
$\G$, then $a$ is necessarily invertible.

\Proof Suppose, by way of contradiction, that $a$ is non-invertible.  By \lcite{\Invertibility} there exists an
irreducible representation $\pi$ of $C^*(\G)$ such that $\pi(a)$ is non-invertible.  Employing \cite[Theorem 2.1]{\IW} we
have that, for some $x \in\Gz$, there exists an irreducible representation $L$ of $\CGx$ such that $\pi$ and
$\sysind{\G(x)}{\G}L$ share null spaces.

Since $\G$ is amenable we have that $\G(x)$ is also amenable by
  \cite[Proposition 5.1.1]{\ClaireRenault},
  and hence that $L$ is weakly contained in the left-regular representation.  We may therefore employ
\lcite{\Comparison} to conclude that
  $$
  \Ker(\lambda_x) \subseteq\Ker\big(\sysind{\G(x)}{\G}L\big) = \Ker(\pi).
  $$
  By hypothesis $a$ is invertible modulo $\Ker(\lambda_x)$, and hence it must also be invertible modulo $\Ker(\pi)$, a
contradiction.  \endProof

\section Strictly norming family of representations

A family $\ScriptCapitalF$ of representations of a C*-algebra $A$ is often called \"{norming}, when
  $$
  \Vert a\Vert = \sup_{\pi \in \ScriptCapitalF} \Vert\pi(a)\Vert \for a\in A.
  \eqmark Norming
  $$
  As an example, the family $\{\lambda_x\}_{x\in\Gz}$ is norming for the reduced groupoid C*-algebra $C_r^*(\G)$, for every
(non-necessarily amenable) groupoid $\G$.  Based on this concept, let us give the following:

\definition A family $\ScriptCapitalF$ of representations of a C*-algebra $A$ will be called \"{strictly norming} when it is norming
and, in addition, the supremum in \lcite{\Norming} is \"{attained} for every $a$ in $A$.

The next result, due to Roch, relates strictly norming and sufficient families in an interesting way.  Its proof is
included for the convenience of the reader and also because it is slightly simpler than the proof given by Roch in
\cite{\Roch}.

\state Theorem \label RochTheorem (\cite[Theorem 5.7]{\Roch}) Let $\ScriptCapitalF$ be a family of non-degenerated representations
of a unital C*-algebra $A$.  Then $\ScriptCapitalF$ is strictly norming if and only if it is sufficient.

\Proof
  Arguing by contradiction, suppose that $\ScriptCapitalF$ is sufficient, but there exists $a \in A$ such that $\Vert\pi(a)\Vert<\Vert a\Vert$, for all $\pi$
in $\ScriptCapitalF$.  Replacing $a$ by $a^*a$, we may assume that $a$ is positive.  For every $\pi$ in $\ScriptCapitalF$, we then have that
  $$
  \def \spc#1{\ \ #1\ \ }
  \sigma\big(\pi(a)\big)\spc \subseteq\big[0,\Vert\pi(a)\Vert\big]\spc \subseteq\big[0,\Vert a\Vert\big).
  $$
  Setting $b = a - \Vert a\Vert$, we then have by the spectral mapping theorem that
  $$
  \def \spc#1{\ \ #1\ \ }
  \sigma\big(\pi(b)\big) \spc = \sigma\big(\pi(a) - \Vert a\Vert\big) \spc = \sigma\big(\pi(a)\big) - \Vert a\Vert\spc \subseteq\big[-\Vert a\Vert,0\big).
  $$
  It follows that $0 \notin\sigma\big(\pi(b)\big)$, and hence that $\pi(b)$ is invertible for every $\pi$ in $\ScriptCapitalF$, but, since $\Vert a\Vert$ belongs to
the spectrum of $a$, we see that $b$ is not invertible, a contradiction.

To verify the ``only if\kern0.8pt" part of the statement, let $a$ be non-invertible.  We thus need to find some $\pi\in\ScriptCapitalF$, such
that $\pi(a)$ is non-invertible.

Since $a$ is non-invertible, then either $a^*a$ or $aa^*$ is non-invertible.  We suppose without loss of generality that
the former is true, that is, that the element $c : = a^*a$ is non-invertible.
  We then have that
  $$
  0 \in\sigma(c) \subseteq\big[0,\Vert c\Vert\big].
  $$
  With $b = \Vert c\Vert-c$, we conclude from the spectral mapping theorem that
  $$
  \def \spc#1{\ \ #1\ \ }
  \Vert c\Vert\spc \in\sigma(b) \spc \subseteq\Vert c\Vert- \big[0,\Vert c\Vert\big] \spc = \big[0,\Vert c\Vert\big],
  $$
  so $\Vert b\Vert= \Vert c\Vert$, and by hypothesis there exists $\pi\in\ScriptCapitalF$, such that $\Vert\pi(b)\Vert= \Vert c\Vert$.  Since $\pi(b)$ is positive, this implies
that $\Vert c\Vert$ lies in its spectrum, which is to say that $\Vert c\Vert-\pi(b)$ is non-invertible, but
  $$
  \Vert c\Vert-\pi(b) = \pi(c),
  $$
  so $\pi(c)$ is non-invertible which implies that $\pi(a)$ is non-invertible.  \endProof

Putting  \lcite{\MainThm} and \lcite{\RochTheorem} together, we therefore deduce the following important consequence:

\state Corollary \label MainCorol Let $\G$ be a second-countable, Hausdorff, {\etale}, amenable groupoid, with $\Gz$ compact.
Then, for every $a\in C^*(\G)$, there exists $x\in\Gz$, such that
  $$
  \Vert a\Vert = \Vert\lambda_x(a)\Vert.
  $$

\vfill

\references

\def \bib #1#2#3#4#5{\bibitem{#1}{#3}{#4}{#5}}

\bib{\Allan}{A} 
  {G. R. Allan}
  {Ideals of vector-valued functions}
  {\it Proc. London Math. Soc. \bf 18 \rm (1968), 193--216}

\bib {\ClaireRenault}{AR}
  {C. Anantharaman-Delaroche and J. Renault}
  {Amenable groupoids}
  {Monographies de L'Enseigne\-ment Mathematique, 36. L'Enseignement
Mathematique, Geneva, 2000. 196 pp}

\bib {\Amelia}{BFK}
  {M. A. Bastos, C. A. Fernandes, Y. I. Karlovich}
  {Spectral measures in C*-algebras of singular integral operators with shifts}
  {J. Funct. Analysis \bf 242 \rm (2007),  86--126}

\bib {\Dougkas}{D}
  {R. G. Douglas}
  {Banach algebra techniques in operator theory}
  {Academic Press, 1972}

\bib {\IW}{IW}
  {M. Ionescu and D. Williams}
  {The generalized Effros-Hahn conjecture for groupoids}
  {\it Indiana Univ. Math. J. \bf 58 \rm (2009), no. 6, 2489--2508}

\bib {\IWTwo}{IW2}
  {M. Ionescu and D. Williams}
  {Irreducible representations of groupoid C*-algebras}
  {\it Proc. Amer. Math. Soc. \bf 137 \rm (2009), 1323--1332}

\bib {\KarlovOne}{K}
  {Y. I. Karlovich}
  {A local-trajectory method and isomorphism theorems for nonlocal C*-algebras}
  {\it Operator Theory: Advances and Applications \bf 170 \rm (2006), 137--166}

\bib {\MRW}{MRW}
  {P. Muhly, J. Renault, and D. Williams}
  {Equivalence and isomorphism for groupoid C*-algebras}
  {\it J. Operator Theory \bf 17 \rm (1987), 3--22}

\bib {\renault}{Re}
  {J. Renault}
  {A groupoid approach to C*-algebras}
  {Lecture Notes in Mathematics vol.~793, Springer, 1980}

\bib {\Roch}{Ro}
  {S. Roch}
  {Algebras of approximation sequences: structure of fractal algebras}
  {Singular integral operators, factorization and applications,
287--310, Oper. Theory Adv. Appl., 142, Birkh\umlaut auser, Basel, 2003.}

\endgroup

\bigskip {\it E-mail address: \tt ruyexel@gmail.com}

\bye